\documentclass{svjour3}
\smartqed
\usepackage{amsmath} 
\usepackage{amsfonts} 
\usepackage[a4paper]{geometry} 
\usepackage{graphicx} 
\usepackage{epstopdf} 
\usepackage{hyperref} 
\usepackage{algorithm,algpseudocode} 
\usepackage{caption,subcaption} 
\usepackage{wrapfig,sidecap} 
\usepackage{comment} 
\usepackage[hang,flushmargin]{footmisc}
\algrenewcomment[1]{\quad\(\triangleright\) #1}
\graphicspath{ {YAF-Figures/FFTfiltering/} }
\algnewcommand\algorithmicinput{\textbf{Input:}}
\algnewcommand\INPUT{\item[\algorithmicinput]}
\algnewcommand\algorithmicoutput{\textbf{Output:}}
\algnewcommand\OUTPUT{\item[\algorithmicoutput]}
\algnewcommand{\LineComment}[1]{\State \(\triangleright\) #1}
\newcommand{\be}{\begin{equation}}
\newcommand{\ee}{\end{equation}}

\usepackage{etoolbox}
\newcommand*{\affaddr}[1]{#1} 
\newcommand*{\affmark}[1][*]{\textsuperscript{#1}}

\begin{document}
\title{Stopping criterion for iterative regularization of large-scale ill-posed problems using the Picard parameter\thanks{This research was supported by the Israel Science Foundation (grant No. 132/14) and by Rosa and Emilio Segr\'{e} Research Award.}}

\titlerunning{Stopping criterion for iterative regularization using the Picard parameter}

\author{Eitan Levin\affmark[1]
\and Alexander Y. Meltzer\affmark[1]}
\authorrunning{E. Levin, A.Y. Meltzer}

\institute{Eitan Levin   \\
\email{eitan.levin@weizmann.ac.il}, \\
\and
Alexander Y. Meltzer \\
\email{alexander.meltzer@weizmann.ac.il}.\\
\affaddr{ \affmark[1] Department of Condensed Matter Physics, Weizmann Institute of Science, 76100 Rehovot, Israel}}

\maketitle
\setcounter{footnote}{1}
\begin{abstract}
We propose a new stopping criterion for Krylov subspace iterative regularization of large-scale ill-posed inverse problems. Our stopping criterion accurately filters the data using a generalization of the Picard parameter that was originally introduced for direct regularization of small-scale problems. In the one dimension we filter the data in the discrete Fourier transform (DFT) basis using the Picard parameter, which separates noise-dominated Fourier coefficients from the signal-dominated ones. For two-dimensional problems we propose a novel vectorization scheme of the Fourier coefficients of the data based on the Kronecker product structure of the two-dimensional DFT matrix, which effectively reduces the problem to one dimension. At each iteration we compute the distance between the data reconstructed from the iterated solution and the filtered data, terminating the iterations once this distance begins to increase or to level off. The accuracy and robustness of the proposed method is demonstrated by several numerical examples and a MATLAB-based implementation is provided\footnote{\hspace{1pt} A MATLAB-based implementation of the described algorithms can be found at : \\ \href{https://www.weizmann.ac.il/condmat/superc/software/}{https://www.weizmann.ac.il/condmat/superc/software/}.\label{CodeAddr}}.
\end{abstract}

\section{Introduction}\label{sec:Intro}
Iterative methods based on projections of the solution onto a Krylov subspace are often used to solve large-scale linear ill-posed problems \cite{RegParamItr,Chung2008,Novati2013,Bazan2010,Hochstenbach2010,projRegGenForm,Gazzola2015}, \cite[Chap. 6]{RankDeff}, \cite[Chap. 6]{HansenInsights}. Such large problems arise in a variety of applications including image deblurring \cite{Nagy2004,Yuan2007,SpecFiltBook,projRegGenForm,Chung2008} and machine learning \cite{Ong2004,Martens2010,Freitas2005,Ide2007}. Krylov methods iteratively project the solution onto a low-dimensional subspace and solve the resulting small-scale problem using standard procedures such as the QR decomposition \cite{LSQR}. These methods are therefore attractive for the solution of large-scale problems that cannot be solved directly as well as for problems perturbed by noise, since the projection onto a Krylov subspace possesses a regularizing effect \cite{Hnetynkova2009}. Accurate solution using iterative procedures requires  stopping the process close to the optimal stopping iteration. In this paper we develop a general stopping criterion for Krylov subspace regularization, which we particularly apply to the Golub-Kahan bidiagonalization (GKB), also frequently referred to as Lanczos bidiagonalization \cite{Golub1965}.

Consider the problem
\be\label{eq:Ax=b} Ax=b,\ee
where the matrix \(A\in\mathbb{R}^{m\times n}\) is large and ill-conditioned and the data vector \(b=b_{true} + n\) constitutes the true data \(b_{true}\) perturbed by an additive white noise vector \(n\). We are interested in approximating the least-squares solution of the problem
\be\label{eq:LSsolution}
\min_x ||b_{true}- Ax||^2
\ee
without knowledge of \(b_{true}\). This can be done by minimizing the projected least-squares (PLS) problem
\be\label{eq:GKBprob}
\min_x ||b-Ax||^2,\quad \text{such that}\quad x\in\mathcal{K}_k(A^TA,A^Tb),
\ee
where
\be\label{eq:KrylovSubspace}
\mathcal{K}_k(A^TA,A^Tb)=\text{span}\{A^Tb,\ldots,\left(A^TA\right)^{k-1}A^Tb\},
\ee
is the Krylov subspace generated using the GKB process for each iteration \(k\).
The regularizing effect of the projection onto \(\mathcal{K}_k(A^TA,A^Tb)\) then dampens the noise in \(b\), allowing us to approximate \(x_{true}\) \cite{RegParamItr,Hnetynkova2009}. The reason for this regularization effect is that in the initial iterations the basis vectors spanning the Krylov subspace are smooth and so is the projected solution. However, for large iteration numbers \(k\) the accuracy of the projected solution decreases as the basis vectors become corrupted by noise. By applying GKB to \ref{eq:Ax=b}, we thus obtain a sequence of iterates \(\{x^{(k)}\}\) whose error \(||x^{(k)}-x_{true}||\) initially decreases with increasing iterations and then goes up sharply. This behavior of GKB is termed semi-convergence \cite[Sect. 6.3]{RankDeff}. It is therefore crucial to develop a reliable stopping rule by which to terminate the iterative solution of \ref{eq:GKBprob} before noise contaminates the solution. For this purpose, a number of stopping criteria have been proposed including the L-curve \cite{RegParamItr,LCurve}, the generalized cross validation (GCV) \cite{RegParamItr}, \cite[Sect. 7.4]{RankDeff} and the discrepancy principle \cite{RegParamItr}. However, both the GCV and the L-curve methods are inaccurate for determination of the stopping iteration in a significant percentage of cases, as has been demonstrated in \cite[Sect. 7.2]{Hansen2006} for the former and below in \ref{sec:NumEx} for the latter. The discrepancy principle, on the other hand, requires \emph{a priori} knowledge of the noise level and is highly sensitive to it \cite[Sect. 4.1.2]{Bauer2011}. Recently, a new stopping criterion called the normalized cumulative periodogram (NCP), which is based on a whiteness measure  of the residual  vector \(r^{(k)}=b-Ax^{(k)}\), was proposed in \cite{Hansen2006}. While this method outperforms the above-mentioned alternatives, we nevertheless show that its results are inconsistent for some of our numerical problems.

Estimation of the stopping iteration using the GCV or the L-curve method requires projection of the solution onto the subspace \(\mathcal{K}_k(A^TA,A^Tb)\), which depends on the noisy data vector \(b\). In contrast to the original problem \ref{eq:Ax=b} where the noise is entirely contained within data vector \(b\), while coefficient matrix \(A\) is noiseless, the coefficient matrix in the projected problem is contaminated by noise. It is thus nontrivial to generalize standard parameter-choice methods to the problem \ref{eq:GKBprob}, as their usage may result in suboptimal solutions due to the fact that they do not account for noise in the coefficient matrix. To overcome this problem, we estimate the optimal stopping iteration for GKB using the Data Filtering (DF) method originally proposed and briefly discussed in \cite{Levin2016}. Using the DF method the stopping iteration is selected as the one for which the distance \(||\widehat b-Ax^{(k)}||\) between the filtered data \(\widehat b\) and the data reconstructed from the iterated solution \(Ax^{(k)}\) is either minimal or levels-off. In \cite{Levin2016} the filtered data \(\widehat b\approx b-n\) is obtained by separating the noise from the data using the so-called Picard parameter \(k_0\), above which the coordinates of the data in the basis of the left singular vectors of \(A\) are dominated by noise. The approximation of the unperturbed data \(\widehat b\) is then obtained by setting the coordinates of \(b\) in that basis to zero for \(k>k_0\).

Computing the SVD of \(A\) for large-scale problems is not feasible and the specific basis in which the authors of \cite{Levin2016} achieve a separation between signal and noise is therefore unavailable. To overcome this problem we propose to replace the SVD basis used in \cite{Levin2016} by the basis of the discrete Fourier transform (DFT) and show that we can achieve a similar separation of signal from noise in one and in two dimensions. It is well known however, that the DFT assumes the signal to be periodic, and applying it to a non-periodic signal results in artifacts in the form of fictitious high-frequency components that cannot be distinguished from the noise in the DFT basis. We prevent these high-frequency components from appearing by using the Periodic Plus Smooth (PPS) decomposition \cite{Moisan2011}, which allows us to write the signal as a sum, \(b=p+s\), of a periodic component \(p\) and a smooth component \(s\). The periodic component is compatible with the periodicity assumption of the DFT and does not produce high-frequency artifacts, while the smooth component does not need to be filtered. We then have to filter only the Fourier coefficients of the periodic component.

For two-dimensional problems, the Fourier coefficients of the data require a vectorization prior to estimation of the Picard parameter. The coefficients are usually vectorized by order of increasing spatial frequency \cite{Hansen2006}. Here we propose an alternative vectorization, ordered by increasing value of the product of the spatial frequencies in each dimension, which enables a more accurate and consistent estimation of the Picard parameter. We demonstrate that such ordering is equivalent to the sorting of a Kronecker product of two vectors of spatial frequencies and stems from a corresponding Kronecker product structure of the two-dimensional DFT matrix. This approach is also analogous to reordering the SVD of a separable blur as discussed in e.g. \cite[Sect. 4.4.2]{SpecFiltBook}. We demonstrate that a filter based on the proposed ordering performs similarly to or outperforms its spatial frequencies-based counterpart in all our numerical examples, allowing termination of the iterative process closer to the optimal iteration. The new filtering procedure is simple and effective, and can be used independently of the iterative inversion algorithm.

Hybrid methods, which replace the semi-convergent PLS problem \ref{eq:GKBprob} with a convergent alternative, have received significant attention in recent years  \cite{RegParamItr,Chung2008,Novati2013,Hochstenbach2010,Bazan2010,Gazzola2015,Chung2015}. These methods combine Tikhonov regularization with a projection onto \(\mathcal{K}_k(A^TA,A^Tb)\), replacing problem \ref{eq:GKBprob} with the projected Tikhonov problem
\be\label{eq:TikhMinProb} \min_x ||b-Ax||^2 + \lambda^2||Lx||^2\quad \text{such that}\quad x\in\mathcal{K}_k(A^TA,A^Tb),\ee
where \(L\) is a regularization matrix and \(\lambda\) is a regularization parameter that controls the smoothness of the solution. In this paper, we follow the authors of \cite{RegParamItr,Chung2008,Novati2013,Hochstenbach2010,Bazan2010,Gazzola2015} by considering only the \(L=I\) case. We would like to note, however, that a generalization to case \(L\neq I\) was developed and discussed in \cite{Hochstenbach2010}. The advantage of hybrid methods is that given an accurate choice of the value \(\lambda\) at each iteration, the error in the solution of \ref{eq:TikhMinProb} stabilizes for large iterations, contrary to the least-squares problem \ref{eq:GKBprob} for which the solution error has a minimum. However, appropriately choosing the regularization parameter \(\lambda\) at each iteration is a difficult task, since the coefficient matrix in the projected problem becomes contaminated by noise, as in the PLS problem. Hence, standard parameter-choice methods such as the GCV cannot be na\"{\i}vely applied to problem \ref{eq:TikhMinProb}. In practice, the GCV indeed fails to stabilize the iterations in a large number of cases as reported in \cite{Chung2008} and \cite[Sect. 5.1.1]{Bazan2010}. To overcome this problem, the authors of \cite{Chung2008} proposed the weighted GCV (W-GCV) method which incorporates an additional free weight parameter, chosen adaptively at each iteration and shown to significantly improve the performance of the method. We demonstrate however, using several numerical examples, that the results of the W-GCV method are still suboptimal.

This paper is organized as follows. In Sect. \ref{sec:DirectReg} we summarize results from the Tikhonov regularization of \ref{eq:Ax=b} that we extend to the PLS problem \ref{eq:GKBprob}. In Sect. \ref{sec:DFTfilter}, we present our filtering technique, based on the Picard parameter in the DFT basis for one- and two-dimensional problems. In Sect. \ref{sec:GKBinvert}, we present the GKB algorithm and formulate our stopping criterion. In this section we also discuss hybrid methods that combine projection with Tikhonov regularization. Finally, in Sect. \ref{sec:NumEx} we give numerical examples which demonstrate the performance of the proposed stopping criterion, and compare it to the L-curve, NCP and W-GCV methods.

\section{Tikhonov regularization}\label{sec:DirectReg}
We begin with a description of our parameter-choice method, detailed in \cite{Levin2016}, for standard direct Tikhonov regularization of \ref{eq:Ax=b} using the Picard parameter, which represent the starting point of our derivation. The Tikhonov regularization method solves the ill-posed problem \ref{eq:Ax=b} by replacing it with the related, well-posed counterpart
\be\label{eq:TikhProb} \min_x ||b-Ax||^2 + \lambda^2||x||^2 \implies \left(A^TA+\lambda^2I\right)x=b.\ee
The solution of \ref{eq:TikhProb} can be written in a convenient form, using the singular value decomposition (SVD) of \(A\), given by
\be\label{eq:SVDofA} A = U\Sigma V^T,\ee
where \(U\in\mathbb{R}^{m\times m}\) and \(V\in\mathbb{R}^{n\times n}\) are orthogonal matrices. For simplicity, let the \(j\)th columns of \(U\) and \(V\) be denoted by \(u_j\) and \(v_j\), respectively, and the \(j\)th singular value of \(A\) by \(\sigma_j\). Furthermore, let \(\beta_j= u_j^Tb\) be the \(j\)th Fourier coefficient of \(b\) with respect to \(\{u_j\}_{j=1}^m\) and let \(\nu_j = u_j^Tn\) be the Fourier coefficients of the noise. Then, the solution of \ref{eq:TikhProb} can be written as
\be\label{eq:TikhSoln}
x(\lambda) = \sum_{j=1}^{m}\frac{\sigma_j}{\sigma_j^2+\lambda^2}\beta_j v_j.
\ee
It can be shown that in order for solution \ref{eq:TikhSoln} to represent a good approximation to the true solution \(x_{true}\) for some \(\lambda\), the problem must satisfy the discrete Picard condition (DPC) \cite{DPC}. The DPC requires that the sequence of Fourier coefficients of the true data \(\{u_j^Tb_{true}\}=\{\beta_j-\nu_j\}\) decays faster than the singular values \(\{\sigma_j\}\) which, by the ill-conditioning of \(A\), decay relatively quickly. Therefore, the DPC implies that \(\beta_j-\nu_j \approx 0\), or equivalently, \(\beta_j\approx \nu_j\), for \(j\geq k_0\) from some index \(k_0\), termed the Picard parameter, on. In other words, the coefficients of \(b\) with indices larger than \(k_0\) are dominated by noise, while coefficients with smaller indices are dominated by the true data.

To estimate \(\lambda\) we can rewrite \ref{eq:LSsolution} as
\be\label{eq:LSsolution2}
\min_\lambda ||b_{true}- Ax(\lambda)||^2,
\ee
but since \(b_{true}\) is not known we suggest to replace it with the filtered field \(b_{true}\approx \hat b\), as in the DF method \cite{Levin2016}. The DF method sets the regularization parameter \(\lambda\) for \ref{eq:TikhProb} to be the minimizer of the distance function
\be
\label{eq:DistNorm}
\min_\lambda ||\widehat{b}-Ax(\lambda)||^2,
\ee
between the data \(Ax(\lambda)\) reconstructed from the solution \ref{eq:TikhProb}, and the filtered data \(\widehat b\). To obtain the filtered data \(\widehat b\) we remove the noise-dominated coefficients from the expansion of \(b\) in basis \(\{u_j\}\) so that
\be\label{eq:PicFiltSVD}
\widehat b = \sum_{j=1}^{k_0-1}\beta_ju_j.
\ee
The Picard parameter \(k_0\), can be found by detection of the levelling-off of the sequence
\be\label{eq:VxDefn}
V(k) = \frac{1}{m-k+1}\sum_{j=k}^m \beta_j^2,
\ee
which is shown to decrease on average until it levels-off at \(V(k_0)\simeq s^2\), where \(s^2\) is the variance of the white noise. The detection is done by setting \(k_0\) to be the smallest index satisfying
\be\label{eq:PicIndCondSVD}
\frac{|V(k+h)-V(k)|}{V(k)} \leq \varepsilon,
\ee
for some step size \(h\) and a bound on the relative change \(\varepsilon\). The above estimation of \(k_0\) is stated to be robust to changes in \(h\) and \(\varepsilon\), with the values \(\varepsilon\in[10^{-3},10^{-1}]\) and \(h\in[\lfloor\frac{m}{100}\rfloor,\lceil\frac{m}{10}\rceil]\) working consistently well \cite{Levin2016}.

Unfortunately, for large-scale problems, computing the SVD of \(A\) given by \ref{eq:SVDofA} is unfeasible in general and therefore the above separation of noise from signal in the SVD basis is unobtainable. In the next section, we propose a similar filtering procedure for \(b\), which utilizes the basis of the DFT instead of SVD. This basis satisfies an analog of the DPC and is effective for large-scale applications.

\section{The DFT filter}\label{sec:DFTfilter}
In this section we replace the SVD basis discussed in Sect. \ref{sec:DirectReg} with the DFT basis. Since computing the Fourier coefficients with respect to the DFT basis can be done efficiently, the proposed procedure remains computationally cheap even for large-scale problems.
We begin by noting that the true data \(b_{true}\) is generally smooth and therefore is dominated by low-frequencies. This is true for cases such as image deblurring problems and problems arising from discretization of integral equations, where the coefficient matrix \(A\) has a smoothing effect and hence, \(b_{true} = Ax_{true}\) is smooth even if \(x_{true}\) is not, see \cite{Hansen2006,Hansen2008} \cite[Sect. 5.6]{SpecFiltBook}. Furthermore, the SVD basis \(\{u_j\}\) is usually similar to that of the DFT,as shown in \cite{Hansen2006}, where the authors demonstrate that vectors \(u_j\) corresponding to small indices \(j\) are well represented by just the first few Fourier modes. In contrast, vectors \(u_j\) with large indices \(j\) are shown to include significant contributions from high frequency Fourier modes. These observations suggest that we can replace the SVD basis with the Fourier basis, so that the role of the decreasing singular values in the ordering of the basis vectors is played by the increasing Fourier frequencies. For our procedure to be valid we expect the DFT coefficients of \(b_{true}\) to satisfy an analog of the DPC and therefore to decay to zero as the frequency increases.

For an image \(B\) of size \(M\times N\) we use the two-dimensional DFT
\be\label{eq:DFT2D}
\text{DFT2}[B] = \mathcal{F}_M^*B\overline{\mathcal{F}}_N,
\ee
where
\be\label{Eq:DFTmat} \left(\mathcal{F}_m\right)_{j,k} = \frac{1}{\sqrt{m}}e^{i2\pi(j-1)(k-1)/m},\ee
is the unitary DFT matrix of size \(m\times m\), \(\overline{X}\) denotes complex conjugation, and \(i=\sqrt{-1}\). Note that \ref{eq:DFT2D} reduces to the one-dimensional DFT if \(N=1\). The data vector \(b\) in \ref{eq:Ax=b} is then obtained by vectorizing the matrix \(B\) by stacking its columns one upon the other so that \(b=\text{vec}(B)\), where \(\text{vec}(\cdot)\) denotes the above vectorization scheme and \(m=MN\) is the resulting length of \(b\). However, the Fourier coefficients found in \ref{eq:DFT2D} cannot be used directly for our purposes because a na\"{\i}ve application of DFT to a non-periodic signal causes artifacts in the frequency domain. Specifically, DFT assumes that the data to be transformed is periodic and therefore application of the DFT to a non-periodic data leads to discontinuities at the boundaries. In the frequency domain, these discontinuities take the form of large high-frequency coefficients \cite{Moisan2011}. Thus Fourier coefficients of smooth but non-periodic true data do not satisfy the DPC as we require. To circumvent this difficulty, we propose to use the Periodic Plus Smooth (PPS) decomposition introduced in \cite{Moisan2011}. The PPS decomposition decomposes an image into a sum
\be\label{eq:PPS} B = P + S,\ee
of a periodic component \(P\) very similar to the original one but that smoothly decays towards its boundaries, and a smooth component \(S\) that is nonzero mainly at the boundaries to compensate for the decaying \(P\).
To compute the PPS of \(B\), we define
\be\label{eq:PerGapImg}
\begin{aligned}
&V_1(j,k) = \left\{\begin{array}{cc} B(M-j+1,k)-B(j,k), & \text{if } j=1 \text{ or } j=M,\\ 0, & \text{otherwise},\end{array}\right.,\\
&V_2(j,k) = \left\{\begin{array}{cc} B(j,N-k+1)-B(j,k), & \text{if } k=1 \text{ or } k=N,\\ 0, & \text{otherwise},\end{array}\right.,
\end{aligned}
\ee
and \(V = V_1 + V_2\). Then, the two-dimensional DFT of the smooth component \(S\) is given by
\be\label{eq:DFTS}
\text{DFT2}[S](j,k) = \left\{\begin{array}{cc} 0, & \text{if } j=k=1,\\ \frac{\text{DFT2}[V](j,k)}{2\cos\left(\frac{2\pi(j-1)}{M}\right)+2\cos\left(\frac{2\pi(k-1)}{N}\right)-4}, & \text{otherwise},\end{array}\right.
\ee
which can be inverted to obtain \(S\) and \(P=B-S\). The subsequent filtering procedure uses only the periodic component \(P\), which contains all the noise as \(S\) is always smooth.

In order to filter the vectorized Fourier coefficients
\be\label{eq:VecFourierCoeffs2D} \beta = \text{vec}\left(\text{DFT2}[P]\right),
\ee
by using the Picard parameter, as in Sect. \ref{sec:DirectReg}, we first have to rearrange \(\beta\) so that the first coefficients correspond to the true data while the last are dominated by noise. In \cite{Hansen2006} the coefficients were arranged in order of increasing spatial frequency. Specifically, the basis of the two dimensional Fourier transform is a plane wave given by
\be\label{eq:2DFTbasis}
f((j,k),(s,l)) = \exp\left\{-i2\pi \left[\frac{(j-1)(k-1)}{M}+\frac{(s-1)(l-1)}{N}\right]\right\} = \exp\left\{-i2\pi \mathbf{k}\cdot\mathbf{x}\right\},
\ee
where \(\mathbf{k}=\left(\frac{j-1}{M},\frac{s-1}{N}\right)^T\) is the frequency vector, \(j=0,1,...,M\), \(s=0,1,...,N\) and \(\mathbf{r}=(k-1,l-1)^T\) is the spatial vector. The components of \ref{eq:VecFourierCoeffs2D} are arranged in order of increasing magnitude of the spatial frequency \(\mathbf{k}\), given by \(|\mathbf{k}|^2=(j-1)^2/M^2+(s-1)^2/N^2\). We refer to this ordering as the elliptic ordering since the contours of the spatial frequency function \(|\mathbf{k}|^2\) are ellipses centered about the zero frequency. However, use of this arrangement in our numerical experiments causes some results to be highly suboptimal.

An alternative ordering of the Fourier coefficients allows us to overcome this problem. Specifically, we utilize the Kronecker product structure of the two-dimensional Fourier transform, which can be written as a matrix-vector multiplication with \(b\) as
\be\label{eq:DFT2kron}
\text{vec}\left(\text{DFT2}[B]\right) = \left(\mathcal{F}^{(2)}_{M,N}\right)^*b.
\ee
Here
\be\label{eq:2DForuierMat}
\mathcal{F}^{(2)}_{M,N} = \mathcal{F}_M\otimes\mathcal{F}_N
\ee
is the 2D Fourier transform matrix and '\(\otimes\)' denotes the Kronecker product defined as
\be\label{eq:KronDefn} A\otimes B = \left(
                   \begin{array}{ccc}
                     a_{1,1}B & \cdots & a_{1,n}B \\
                     \vdots & \ddots & \vdots \\
                     a_{m,1}B & \cdots & a_{m,n}B \\
                   \end{array}
                 \right).\ee
In view of \ref{eq:2DForuierMat}, we suggest to reorder the Fourier coefficients \(\beta\) in \ref{eq:VecFourierCoeffs2D} according to the ordering permutation \(\pi\) which arranges the Kronecker product
\be\label{eq:freqVec}
\mathbf{f}^{(2)}_{M,N} = \mathbf{f}_M\otimes \mathbf{f}_N\in\mathbb{R}^m,
\ee
in increasing order, where \(\mathbf{f}_M\in\mathbb{R}^M\) and \(\mathbf{f}_N\in\mathbb{R}^N\) are the vectors representing the ordered absolute frequencies of \(\mathcal{F}_M\) and \(\mathcal{F}_N\) respectively. Note that since the frequencies in the two-dimensional Fourier transform are shifted so that the zero frequency is located at the corner of the image, the vectors \(\mathbf{f}_N\) and \(\mathbf{f}_M\) in \ref{eq:freqVec} also need to be correspondingly shifted.

\begin{figure}[!p]
    \centering
    \hspace*{\fill}
    \begin{subfigure}[b]{0.45\textwidth}
        \centering
        \includegraphics[width=\textwidth]{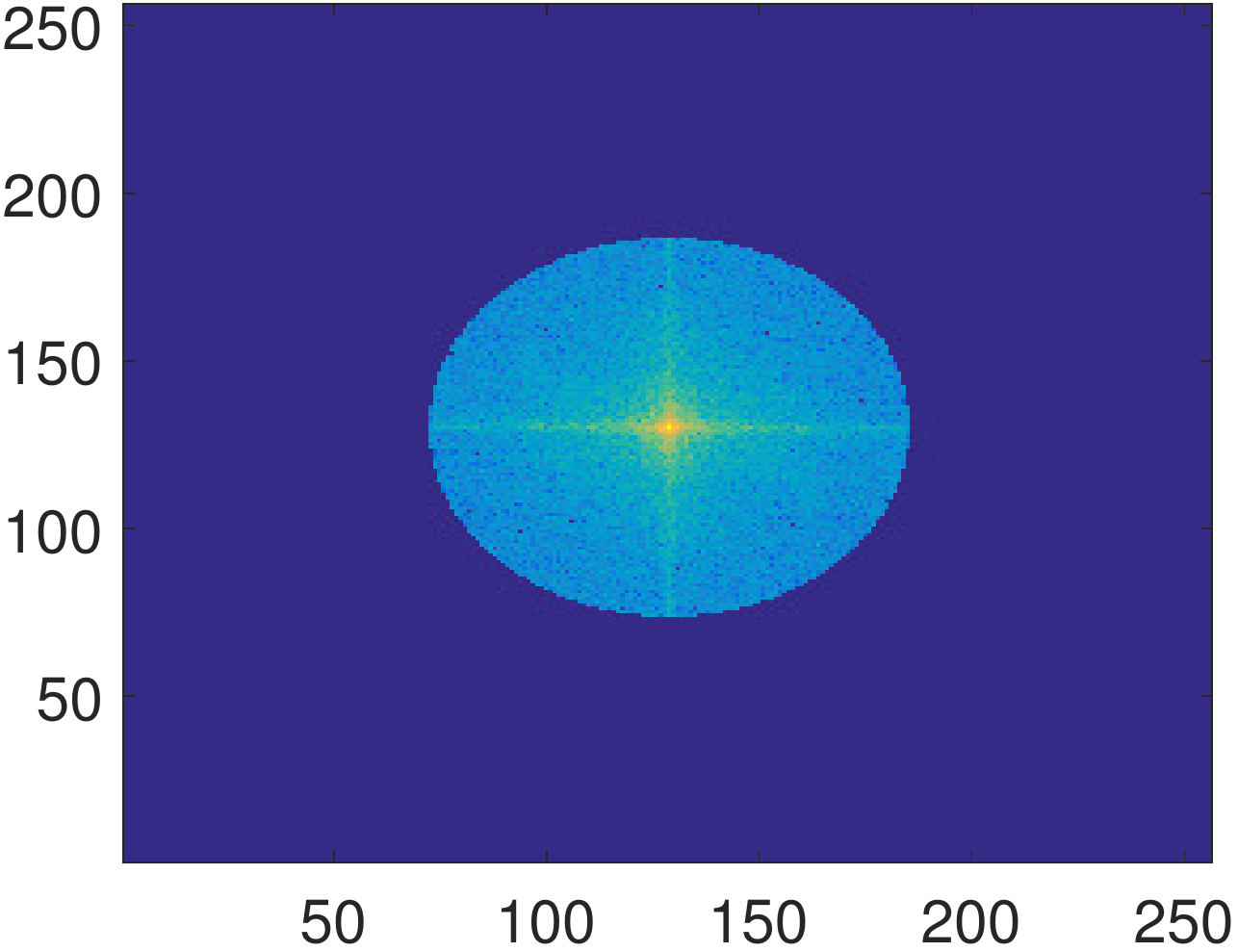}
    \end{subfigure}
    \hfill
    \begin{subfigure}[b]{0.45\textwidth}
        \centering
        \includegraphics[width=\textwidth]{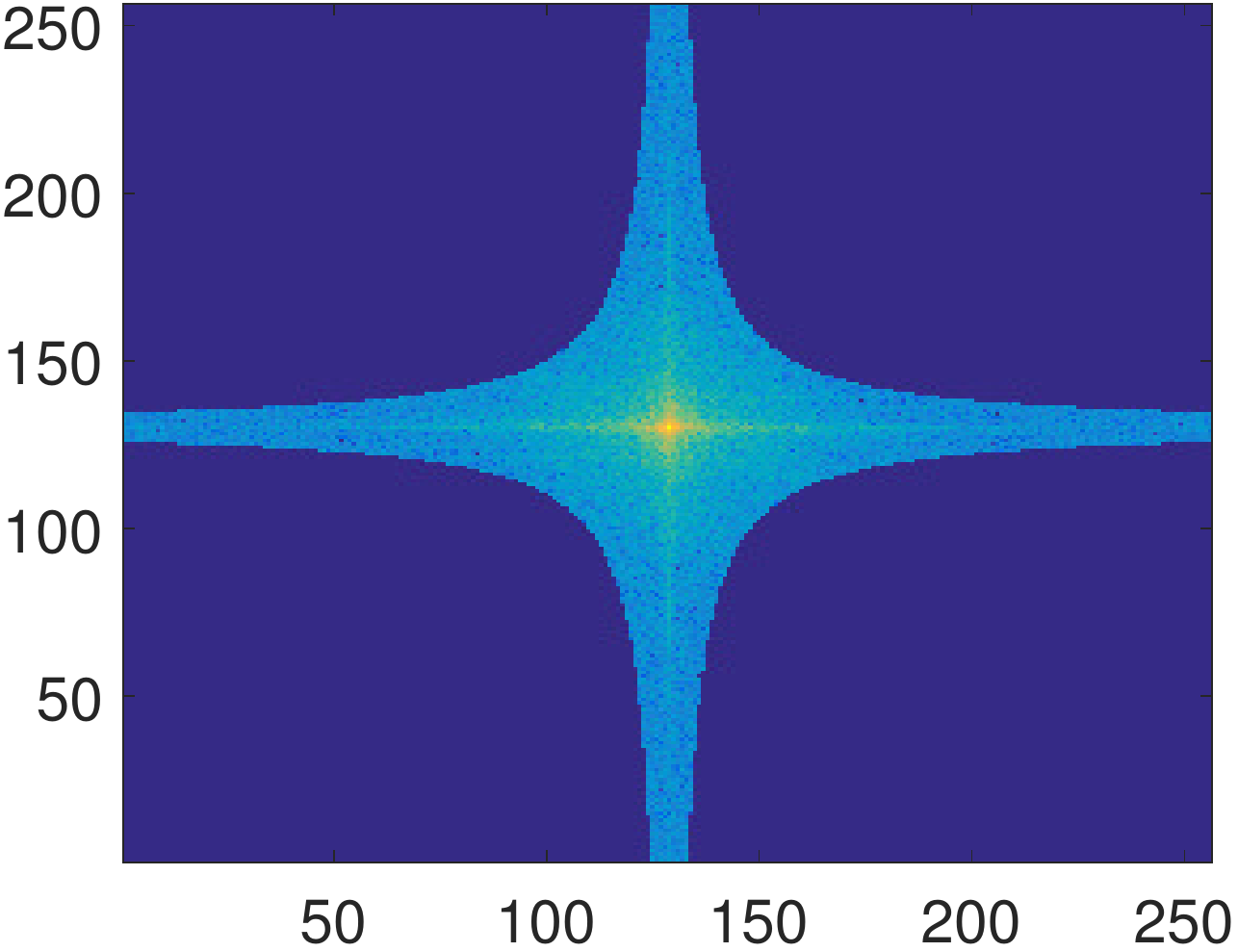}
    \end{subfigure}
    \hspace*{\fill}
\caption{Effect of the filtering procedure introduced in Sect. \ref{sec:DFTfilter} on an image of size \(256\times256\) by the elliptic ordering of \cite{Hansen2006} and the hyperbolic ordering \ref{eq:freqVec} in Fourier domain. The Picard parameter for both methods is \(k_0=10^4\). The zero frequency component is placed at the center of the image.}\label{fig:Masks}
\end{figure}

The vector \(\mathbf{f}^{(2)}_{M,N}\) whose components are products of the absolute values of spatial frequencies is not ordered and in a component-wise form \ref{eq:freqVec} is given by
\be\label{eq:freqVecCompWise}
\left(\mathbf{f}^{(2)}_{M,N}\right)_{M(j-1)+s}=\frac{1}{m}\left\{\begin{array}{ll}
(j-1)(s-1), & 1\leq j \leq \lfloor\frac{N}{2}\rfloor,\ 1\leq s\leq \lfloor\frac{M}{2}\rfloor,\\
(N-j+1)(s-1), & \lfloor\frac{N}{2}\rfloor < j \leq N,\ 1\leq s\leq \lfloor\frac{M}{2}\rfloor,\\
(j-1)(M-s+1), & 1 \leq j \leq \lfloor\frac{N}{2}\rfloor,\ \lfloor\frac{M}{2}\rfloor < s\leq M,\\
(N-j+1)(M-s+1), & \lfloor\frac{N}{2}\rfloor < j \leq N,\ \lfloor\frac{M}{2}\rfloor < s\leq M, \end{array}\right.
\ee
where as above \(m=MN\). We then construct the permutation \(\pi\) such that \(\mathbf{f}^{(2)}_{M,N}(\pi(1:m))\) appears in increasing order and rearrange the coefficients \ref{eq:VecFourierCoeffs2D} to obtain \(\beta\mapsto \beta(\pi(1:m))\). We term this ordering the hyperbolic ordering since the contours of the function \ref{eq:freqVecCompWise} are hyperbolas centered about zero frequency (see \ref{fig:Masks}).

The above rearrangement using the ordering permutation \(\pi\) is analogous to the rearrangement of the SVD decomposition of a separable blur
\be\label{eq:SepBlur} A = A_1\otimes A_2.\ee
where \(A_1\in\mathbb{R}^{N\times N}\), \(A_2\in\mathbb{R}^{M\times M}\) and \(A\in\mathbb{R}^{m\times m}\) \cite[Sect. 2]{Hansen2008}. Letting \(A_1 = U_1\Sigma_1V_1^T\) and \(A_2 = U_2\Sigma_2V_2^T\) be the SVD of \(A_1\) and \(A_2\), the SVD of \(A\) \ref{eq:SVDofA} can be written as
\be
A = \left(\underbrace{U_1\otimes U_2}_{= U}\right)\left(\underbrace{\Sigma_1\otimes \Sigma_2}_{= \Sigma}\right)\left(\underbrace{V_1\otimes V_2}_{= V}\right)^T.
\ee
As in the case of the two-dimensional Fourier transform \ref{eq:2DForuierMat}, even though the singular values of \(A_1\) and \(A_2\) are ordered, those of \(A\) are not \cite[Sect. 4.4.2]{SpecFiltBook}, \cite{Hansen2008}. To be able to use the filter described in Sect. \ref{sec:DirectReg} we must reorder the entries of \(U\), \(\Sigma\) and \(V\) according to decreasing singular values using the ordering permutation \(\pi\) as in \ref{eq:freqVec}.

Once the Fourier coefficients are rearranged we proceed according to the procedure in Sect.  \ref{sec:DirectReg}. Specifically, we form the sequence \(\{V(k)\}_{k=1}^m\) defined in \ref{eq:VxDefn}, estimate the Picard parameter using \ref{eq:PicIndCondSVD} and set to zero the Fourier coefficients with indices larger than \(k_0\) to form the vector \(\widehat\beta(\pi(1:m)) = (\beta(\pi(1)),\ldots,\beta(\pi(k_0-1)),\underbrace{0,\ldots,0}_{m-k_0+1})^T\). We then invert \ref{eq:VecFourierCoeffs2D} using \(\widehat\beta\) instead of \(\beta\) to obtain the filtered periodic component \(\widehat P\) and the filtered image \(\widehat B = \widehat P + S\). The filtered data vector is then obtained as \(\widehat b = \text{vec}(\widehat B)\).

Note that dropping the last coefficients of the data using the two orderings discussed above can also be interpreted as applying one-parameter windows of different shapes in the Fourier domain. Specifically, the elliptic ordering of \cite{Hansen2006} corresponds to setting to zero the Fourier coefficients outside of an ellipse, whereas the hyperbolic ordering corresponds to doing the same outside a hyperbola. This is illustrated in \ref{fig:Masks} where we plot the Fourier transform coefficients of a \(256\times 256\) image with \(k_0=10^4\) for each of the orderings. Viewed from this perspective, the proposed filtering algorithm simply applies a window depending on the parameter \(k_0\) to the DFT of the perturbed image. The difference between our approach and the approach of \cite{Hansen2006} is in the chosen shape of the window.

\section{Iterative inversion using the GKB}\label{sec:GKBinvert}
In this section, we consider the solution of the ill-posed problem \ref{eq:Ax=b} using the GKB algorithm \cite{Golub1965}. This algorithm approximates the subspace, spanned by the first \(k\) largest right singular vectors of \(A\) with the first \(k\) basis vectors of the Krylov subspace \(\mathcal{K}_k(A^TA,A^Tb)\) (see \cite[sect. 6.3.2]{RankDeff} and \cite[sect. 6.3.1]{HansenInsights}). After \(k\) iterations (\(1\leq k\leq n\)), the GKB algorithm yields two matrices with orthonormal columns \(W_k\in\mathbb{R}^{n\times k}\) and \(Z_{k+1}\in\mathbb{R}^{m\times(k+1)}\), and a lower bidiagonal matrix \(B_k\in\mathbb{R}^{(k+1)\times k}\) with the structure
\be\label{eq:BkDefn} B_k = \left(
     \begin{array}{cccc}
       \varrho_1 &   &   &   \\
       \theta_2 & \varrho_2 &   &   \\
         & \theta_3 & \ddots &   \\
         &   & \ddots & \varrho_k \\
         &   &   & \theta_{k+1} \\
     \end{array}
   \right),\ee
such that
\be
\label{eq:LanczosRels}\begin{aligned} &AW_k = Z_{k+1}B_k,\\ &A^TZ_{k+1} = W_kB_k^T + \varrho_{k+1}w_{k+1}e_{k+1}^T,\\ &Z_{k+1}\theta_1e_1 = b.
\end{aligned}
\ee
The GKB algorithm is summarized in \ref{alg:GKB}. We perform a reorthogonalization at each step of the algorithm to ensure that the columns of \(W_k\) and \(Z_{k+1}\) remain orthogonal.
\begin{algorithm}
\caption{Golub-Kahan Bidiagonalization (GKB)}\label{alg:GKB}
\begin{algorithmic}
\INPUT{\(A,b,k\)}\Comment{Coefficient matrix \(A\), data vector \(b\) and number of iterations \(k\)}
\OUTPUT{\(W_k,Z_{k+1},B_k\)}
\LineComment{Initialization:}
\State \(\theta_1 \gets ||b||, \quad z_1\gets b/\theta_1\)
\State \(\varrho_1 \gets ||A^Tz_1||, \quad w_1\gets A^Tz_1/\varrho_1\)
\For{\(j=1,2,\ldots,k\)}
    \State \(p_j \gets Aw_{j}-\varrho_jz_j\)
    \State \(p_j\gets \left(I-Z_jZ_j^T\right)p_j\) \Comment{Reorthogonalization step}
    \State \(\theta_{j+1} = ||p_j||, \quad z_{j+1}\gets p_j/\theta_{j+1}\)
    \State \(q_j \gets A^Tz_{j+1}-\theta_{j+1}w_j\)
    \State \(q_j\gets \left(I-W_jW_j^T\right)q_j\) \Comment{Reorthogonalization step}
    \State \(\varrho_{j+1} = ||q_j||, \quad w_{j+1}\gets q_j/\varrho_{j+1}\)
    \LineComment{Update output matrices:}
    \State \(W_j \gets \left(W_{j-1},\ w_j\right), \quad Z_{j+1} \gets \left(Z_{j},\ z_{j+1}\right)\)
    \State \(B_j \gets \left(\begin{array}{cc} \left(\begin{array}{c} B_{j-1} \\ 0 \end{array}\right), & \left(\begin{array}{c} \mathbf{0} \\ \varrho_j \\ \theta_{j+1} \end{array}\right)\end{array}\right)\)
\EndFor
\end{algorithmic}
\end{algorithm}
It can be shown that the columns of \(W_k\) span the Krylov subspace \ref{eq:KrylovSubspace} and that we can achieve a regularizing effect by projecting the solution onto this subspace \cite[Sect. 2.1]{Gazzola2015}.

Choosing a solution from the column space of \(W_k\), such that \(x=W_ky\), and using the relations \ref{eq:LanczosRels}, we can rewrite the residual norm in \ref{eq:GKBprob} as
\be\rho(k) = ||b-Ax^{(k)}||^2 = ||U_{k+1}\left(\theta_1e_1-B_ky\right)||^2,\ee
and, since \(U_{k+1}\) has orthonormal columns,
\be\label{eq:rewriteLS}
\rho(k) = ||\theta_1e_1-B_ky||^2.
\ee
Hence, solving the PLS problem \ref{eq:GKBprob} amounts to minimizing \ref{eq:rewriteLS}, which is small-scale and can be solved using standard procedures such as the QR decomposition of \(B_k\) as in the LSQR algorithm \cite{LSQR}.

In the absence of noise, the GKB algorithm is terminated once \(\varrho_{k}=0\) or \(\theta_{k+1}=0\).
However, the solution of the PLS problem \ref{eq:rewriteLS} exhibits semi-convergence, whereby the error in the iterates \(||x_{true}-x^{(k)}||\) first decreases as \(k\) increases and then increases sharply well before the above condition is met. This is due to the fact that the columns of \(W_k\) contain increasing levels of noise, as described in \cite{Hnetynkova2009}. Hence, at early iterations the columns are almost noiseless and the solution gets closer to the true one, while at later iterations the solution becomes contaminated by noise. It is thus crucial to appropriately terminate the iterations before the noise becomes dominant.

Usually, standard methods like the GCV or the L-curve are used for estimation of the optimal stopping iteration. However the GCV method assumes that the noise is additive and is fully contained in the data vector \(b\), while the coefficient matrix \(A\) is noiseless. That is indeed the case in the original, large-scale problem \ref{eq:Ax=b} but not in the PLS problem of minimizing \ref{eq:rewriteLS}. In the PLS problem the projected data vector \(Z_{k+1}^Tb=\theta_1e_1\) is a noiseless scaled standard basis vector. The projected coefficient matrix \(B_k = Z_{k+1}^TAW_k\) however, is generated by the Algorithm \ref{alg:GKB} from the noisy data vector \(b\) and depends on the columns of \(W_k\) and \(Z_{k+1}\). Therefore, the derivation of the standard form of the GCV function and the proof of its optimality, such as the ones given in \cite[Thm. 1]{GCV2} no longer apply. The justification for the L-curve method seems to hold for the projected problem, but, as we show in the numerical examples section, for many cases it is far from optimal.

In this paper we propose to use the DF method developed in \cite{Levin2016} to stop the iterative process. The DF method uses the distance between the filtered data \(\widehat b\) and the data reconstructed from the \(k\)th iterate \(Ax^{(k)}\) to characterize the quality of the iterated solution \(x^{(k)}\). Writing the distance as
\be\label{eq:DistNormGKB}
\widehat f(k) = ||\widehat b - Ax^{(k)}||^2,
\ee
and using the methods of Sect. \ref{sec:DFTfilter} to obtain the filtered data \(\widehat b\) we expect  \(\widehat f(k)\) to have a global minimum at or near the optimal iteration. However, \(\widehat f(k)\) may also have local minima, and so we must continue the iterations beyond a potential minimum of \ref{eq:DistNormGKB} to ensure that the function \(\widehat f(k)\) continues to increase. In addition, for problems with very small noise levels, the filter of Sect. \ref{sec:DFTfilter} may not change the data vector \(b\) sufficiently, in which case \(\widehat f(k)\) will not have a minimum. Instead, \(\widehat f(k)\) will flatten after the optimal iteration, since the solution will not become significantly contaminated. To account for all the above behaviors, we propose to terminate the iterations once the relative change in \ref{eq:DistNormGKB} is small enough
\be\label{eq:GKBcond}
\frac{\widehat f(k)-\widehat f(k+1)}{\widehat f(k)} \leq \delta,\quad \text{for \(p\) consecutive iterations}
\ee
for some small bound \(\delta>0\). We emphasize that the numerator in \ref{eq:GKBcond} is not an absolute value and becomes negative if \(\widehat f(k)\) starts to increase. Since we choose \(\delta>0\), the condition \ref{eq:GKBcond} is automatically satisfied if \(\widehat f(k)\) has a minimum. Otherwise the iterations are terminated once \(\widehat f(k)\) becomes very flat. From our experience, the algorithm works well for \(\delta\in [10^{-2}, 10^{-4}]\) and \(p \geq 5\). Finally, we observe that since a discrete iteration number (\(k\in\{1,\ldots,n\}\)) acts as a regularization parameter, the ability of the reconstructed data \(Ax^{(k)}\) in \ref{eq:DistNormGKB} to get closer to the filtered data is limited, making the minimum of \ref{eq:DistNormGKB} less sensitive to the quality of the filtered data \(\widehat b\). Furthermore, the discrete regularization parameter does not limit the accuracy of the solution in practice. As we demonstrate in the numerical examples below, the additional regularization of hybrid methods is unnecessary and only slightly improves the solution, if at all.

\subsection{Regularizing the projected problem}\label{sec:TikhRegProj}
In this section, we briefly discuss hybrid methods for the solution of \ref{eq:Ax=b}, which combine projection onto a Krylov subspace with Tikhonov regularization and are the main competitors to our approach. Instead of attempting to terminate the GKB process at an optimal iteration, hybrid methods replace the PLS problem \ref{eq:GKBprob} with the Tikhonov minimization problem \ref{eq:TikhMinProb}. Using the relations \ref{eq:LanczosRels} and \(x=W_ky\), as in Sect. \ref{sec:GKBinvert}, problem \ref{eq:TikhMinProb} can be rewritten as
\be\label{eq:MinTikhGKB} \min_y ||\theta_1e_1-B_ky||^2 + \lambda^2 ||y||^2,\ee
or as the normal equation
\be\label{Eq:TikhNormalEq0}
\left(B_k^TB_k + \lambda^2 I\right)y = B_k^T\theta_1e_1
\ee
which has the solution
\be\label{Eq:TikhNormalEq}
y_\lambda = \left(B_k^TB_k + \lambda^2 I\right)^{-1}B_k^T\theta_1e_1.
\ee
The solution \ref{Eq:TikhNormalEq} can be rewritten, similarly to Sect. \ref{sec:DirectReg}, using the SVD of \(B_k\) as
\be\label{eq:SVDbk}
B_k = U_k\Sigma_k V_k^T,
\ee
where \(U_k\in \mathbb{R}^{(k+1)\times (k+1)}\) and \(V_k\in\mathbb{R}^{k\times k}\) are orthogonal and \(\Sigma_k\in\mathbb{R}^{(k+1)\times k}\) has the structure
\be\label{eq:SigmaStruct}
\Sigma_k = \left(\begin{array}{c} \text{diag}\{\sigma_1,\ldots,\sigma_k\} \\ \mathbf{0}^T \end{array}\right),
\ee
with the singular values \(\{\sigma_j\}\) arranged in decreasing order.
 Due to the structure of \(B_k\) in \ref{eq:BkDefn} and the fact that \(\varrho_j,\; \theta_j\; >0\) for all relevant iterations, we have \(\text{rank}(B_k)=k\) and \(\sigma_j >0\) for all \(j\leq k\). Denoting the \(j\)th columns of \(U_k\) and \(V_k\) as \(u^{(k)}_j\) and \(v^{(k)}_j\) respectively, the Tikhonov solution \ref{Eq:TikhNormalEq} can be written similarly to \ref{eq:TikhSoln} as
\be\label{eq:TikhSolnGKB} y_{\lambda} = \theta_1\sum_{j=1}^k \frac{\sigma^{(k)}_ju_j^{(k)}(1)}{\left(\sigma^{(k)}_j\right)^2+\lambda^2}v^{(k)}_j,\ee
where \(u^{(k)}_j(1)=e_1^Tu^{(k)}_j\).

By appropriately choosing \(\lambda\) at each iteration we can, in theory, filter out noise added to the least-squares solution in \ref{eq:GKBprob} at higher iterations and thus, stabilize the error and make the final solution independent of the stopping iteration. Regularization in hybrid methods may also filter out noise that is not filtered by projection alone. Nevertheless we argue that this additional filtering has a negligible effect and that choosing \(\lambda\) appropriately at each iteration presents a significant challenge. Specifically, determination of \(\lambda\) for hybrid methods is usually done using standard procedures originally developed for direct regularization, the most popular of which is the GCV \cite{splines,GCV2}. These procedures assume a noiseless coefficient matrix \(A\) \ref{eq:MinTikhGKB}, whereas the hybrid methods project the solution into a noise-dependent Krylov space, thus also contaminating the projected coefficient matrix, similarly to the PLS problem. Therefore, the GCV method is not expected to produce optimal solutions for hybrid methods, as was indeed demonstrated in \cite[Sect. 4]{Chung2008}. The W-GCV method proposed in \cite{Chung2008} attempts to overcome the above shortcoming of the GCV by introducing an additional free parameter and choosing it adaptively at each iteration. However, as we show in the numerical examples in Sect. \ref{sec:NumEx}, the W-GCV method still produces suboptimal results in many cases.

In all of our numerical examples we observe that the minimum errors achievable using PLS and hybrid methods are almost identical and therefore any additional filtering achievable by hybrid methods is negligible. To explain this we note that at early iterations, the vectors spanning the Krylov subspace \ref{eq:KrylovSubspace}, and hence also the solutions projected onto them, are typically very smooth and do not contain noise \cite{Hnetynkova2009}. Therefore, little to no regularization is required at this stage and we expect \(\lambda\approx 0\), making the Tikhonov problem \ref{eq:TikhMinProb} equivalent to the least-squares problem \ref{eq:GKBprob}. Only after the basis vectors spanning \ref{eq:KrylovSubspace} become contaminated by noise does the solution require regularization. At this stage the optimal regularization parameter \(\lambda\) increases to some non-negligible, noise-dependent value that keeps the error of the regularized solution approximately constant, while the error of the unregularized PLS solution increases sharply. We demonstrate in the following section, using several numerical examples, that the optimal solution of the hybrid method occurs while the noise that contaminates it is very small and so in most typical cases the optimal solutions of \ref{eq:GKBprob} and \ref{eq:TikhMinProb} are very close to each other. This was also demonstrated in the numerical examples in  \cite[Sect. 4]{Chung2015} and \cite{Chung2008}. Therefore, in most practical problems the only significant advantage of hybrid methods over simple GKB stopping criteria is in the ability to stabilize the iterations, making them less sensitive to the stopping iteration. This also implies that having a reliable stopping criterion for PLS obviates the need for a hybrid method in these cases.

\section{Numerical examples}\label{sec:NumEx}
In this section we demonstrate the performance of the proposed method using seven test problems from the \texttt{Matlab} toolbox \texttt{RestoreTools} \cite{Nagy2004}: \texttt{satellite}, \texttt{GaussianBlur440}, \texttt{AtmosphericBlur50}, \texttt{Grain}, \texttt{Text}, \texttt{Text2} and \texttt{VariantMotionBlur\_large}. Each of these problems includes a different blur, \(A\), and a different image, \(x_{true}\), to reconstruct. To generate the data vector \(b\), we form the true data, \(b_{true} = Ax_{true}\) and perturb it with white Gaussian noise of variance \(s^2=\alpha\max|b_{true}|^2\) where the noise level \(\alpha\) takes on three values for each problem \(\alpha\in\{10^{-2},10^{-4},10^{-6}\}\). We apply the inversion procedure to each test problem 100 times, each time using a different noise realization.

To implement our stopping criterion, we set \(h=\lceil\frac{m}{100}\rceil\) and \(\varepsilon=10^{-2}\) in \ref{eq:PicIndCondSVD} and also \(p=5\) and \(\delta=2\times10^{-3}\) in \ref{eq:GKBcond}. As mentioned above, however, our numerical results are robust and wouldn't change much for a wide range of \(h\), \(\varepsilon\), \(p\) and \(\delta\). In the numerical tests we compare our stopping criterion with the L-curve criterion \cite{RegParamItr,LCurve,Calvetti2000}, \cite[Chap. 7]{RankDeff}, the NCP method for the PLS problem \cite{Hansen2006,Rust2008}, and the hybrid W-GCV method \cite{Chung2008}. The L-curve method constitutes finding the point of maximum curvature on the so-called L-curve, defined as \((||r^{(k)}||,\; ||x^{(k)}||)\), where \(r^{(k)}=b-Ax^{(k)}\) is the residual vector. To do so we use the function \verb!l_corner! from Hansen's \texttt{Regularization Tools} toolbox \cite{RegTools}. Using the L-curve method, we terminate the iterations once the chosen iteration number either stays the same or decreases for \(p=5\) consecutive iterations, signifying that the corner of the L-curve is found.

The NCP method is based on calculating a whiteness measure of the residual vector \(r^{(k)}\) at each iteration \cite{Hansen2006}. The stopping iteration is chosen as the one at which the residual vector \(r^{(k)}\) most resembles white noise, as follows. The vector \(r^{(k)}\) is reshaped into an \(M\times N\) matrix \(R^{(k)}\) satisfying \(r^{(k)} = \text{vec}\left(R^{(k)}\right)\), and the quantity \(\widehat R = |\text{DFT2}(R^{(k)})|\) is defined to be the absolute value of its two-dimensional Fourier transform. Since \(R^{(k)}\) is a real valued signal, it follows that \(\widehat R\) is symmetric about \(q_1 = \lfloor M/2 \rfloor + 1\) and \(q_2 = \lfloor N/2 \rfloor +1\), so that \(\widehat R_{j,k} = \widehat R_{M-j,k}\) for \(2\leq j\leq q_1-1\) and \(\widehat R_{j,k} = \widehat R_{j,N-k}\) for \(2\leq k\leq q_2-1\). Consequently, only the first quarter of \(\widehat R\), which can be written as \(\widehat T = \widehat R(1:q_1,1:q_2)\) using \texttt{Matlab} notation, is required for the analysis. Vector \(\widehat t\) is then obtained by vectorizing \(\widehat T\) using the elliptical parametrization defined in Sect. \ref{sec:DFTfilter}. The NCP of \(R^{(k)}\) is defined as the vector \(c(R^{(k)})\) of length \(q_1q_2-1\) with components\footnote{In \cite{Hansen2006}, the authors assume the problem is square so that \(M=N\) and \(q_1=q_2\).}
\be\label{eq:NCPdefn}
c(R^{(k)})_j = \frac{||\widehat t(2:j+1)||_1}{||\widehat t(2:q_1q_2)||_1},
\ee
where the dc component of \(R^{(k)}\) is not included in the NCP. We note that it is argued in \cite{Rust2008} that the NCP should include the dc component of \(R^{(k)}\) since it captures the deviation from zero mean white noise. However, we found no difference between the two definitions in practice and therefore we shall follow \cite{Hansen2006} and disregard the dc component as in \ref{eq:NCPdefn}. It is shown in \cite{Hansen2006,Rust2008} that for white noise, the NCP should be a straight line from 0 to 1 represented by the vector with components \(s_j = j/(q_1q_2-1)\). Therefore, the whiteness measure is defined as the distance
\be
\label{eq:NCPfunc}
N(k) = ||s-c(R^{(k)})||_1.
\ee
The iterations are terminated once \ref{eq:NCPfunc} reaches its global minimum, signifying that the residual vector at the chosen iteration is the closest to white noise. To implement this method, we compute function \ref{eq:NCPfunc} at each iteration and terminate the iterations once the norm \ref{eq:NCPfunc} increases for \(p=5\) consecutive iterations, just as we do with our own method. We then choose the solution corresponding to the global minimum of the computed \(N(k)\).

In contrast to the above methods, the W-GCV is a hybrid method and solves the Tikhonov problem \ref{eq:TikhMinProb}. It is based on introducing a free parameter to the GCV criterion, as discussed in Sect. \ref{sec:TikhRegProj} and in \cite{Chung2008}. To implement the W-GCV, we use the \texttt{HyBR\_modified} routine provided in the \texttt{RestoreTools} package \cite{Nagy2004} as \texttt{x = HyBR\_modified(A,B,[],HyBRset('Reorth','on'),1)}. Note that we use the reorthogonalization option of the \verb!HyBR_modified! routine to make a fair comparison to our \ref{alg:GKB} that employs full reorthogonalization.

\begin{figure}[!p]
    \centering
    \includegraphics[width=\textwidth]{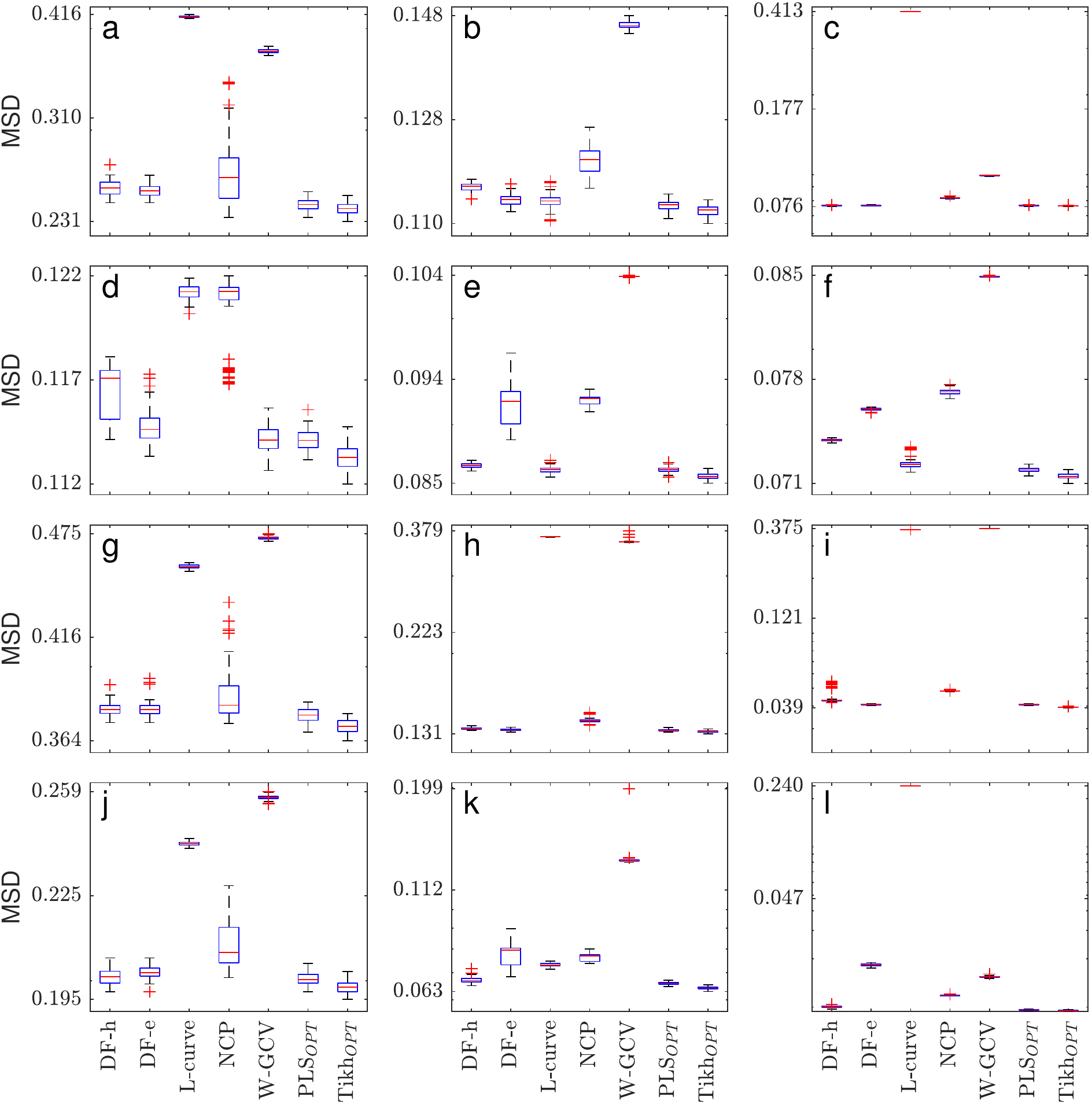}
    \caption{Boxplots of the MSD values obtained for the PLS problem \ref{eq:GKBprob} and the Tikhonov regularized problem \ref{eq:TikhMinProb} with the methods: (1) DF with hyperbolic ordering (DF-h), (2) DF with elliptic ordering (DF-e), (3) L-curve, (4) NCP, (5) W-GCV, (6) minimum MSD for PLS problem (PLS\(_{OPT}\)), (7) minimum MSD solution for Tikhonov problem (Tikh\(_{OPT}\)). The problems presented are \emph{First row}: \texttt{satellite}; \emph{Second row}: \texttt{GaussianBlur440}; \emph{Third row}: \texttt{AtmosphericBlur50}; \emph{Fourth row}: \texttt{Grain}. The noise levels are \emph{First column}: \(\alpha=10^{-2}\); \emph{Second column}: \(\alpha=10^{-4}\); \emph{Third column}: \(\alpha=10^{-6}\).}\label{fig:Boxplots_1-3}
\end{figure}

\begin{figure}[!p]
    \centering
    \includegraphics[width=\textwidth]{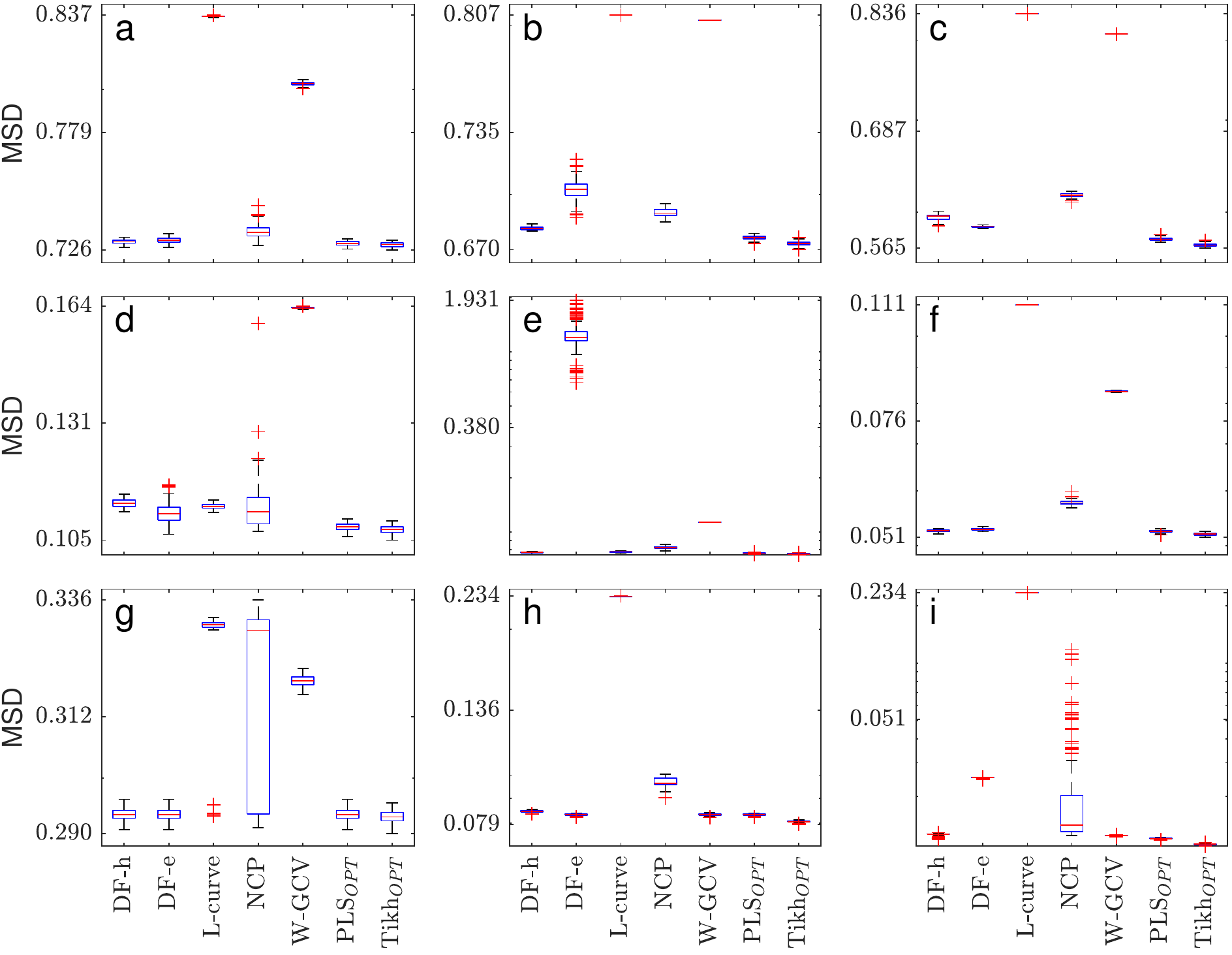}
    \caption{Boxplots of the MSD values obtained for the PLS problem \ref{eq:GKBprob} and the Tikhonov regularized problem \ref{eq:TikhMinProb} with the methods: (1) DF with hyperbolic ordering (DF-h), (2) DF with elliptic ordering (DF-e), (3) L-curve, (4) NCP, (5) W-GCV, (6) minimum MSD for PLS problem (PLS\(_{OPT}\)), (7) minimum MSD solution for Tikhonov problem (Tikh\(_{OPT}\)). The problems presented are \emph{First row}: \texttt{Text}; \emph{Second row}: \texttt{Text2}; \emph{Third row}: \texttt{VariantMotionBlur\_large}. The noise levels are \emph{First column}: \(\alpha=10^{-2}\); \emph{Second column}: \(\alpha=10^{-4}\); \emph{Third column}: \(\alpha=10^{-6}\).}\label{fig:Boxplots_4-6}
\end{figure}

We measure the quality of a solution by computing its Mean-Square Deviation (MSD), defined as
\be\label{eq:MSD} \text{MSD} = \frac{||x_{true} - x||^2}{||x_{true}||^2},\ee
where \(x\) is a solution. We then define the optimal solutions to the PLS problem \ref{eq:GKBprob} and to the hybrid problem \ref{eq:TikhMinProb} as the ones minimizing the MSD to each problem.
We present the results of our simulations as boxplots of the resulting MSD values in \ref{fig:Boxplots_1-3} and \ref{fig:Boxplots_4-6}. The boxplots divide the data into quartiles with the boxes spanning the middle 50\% of the data, called the interquartile range and the vertical lines extending from the boxes span 150\% of the interquartile range above and below it. Anything outside this interval is considered an outlier and is marked with a '+'. Each box also contains a horizontal line marking the median of the data.

Based on the results presented in \ref{fig:Boxplots_1-3} and \ref{fig:Boxplots_4-6} we can make the following observations:
\begin{enumerate}
    \item The hyperbolic ordering with the DF method performs similarly to or better than the corresponding elliptic ordering in all examples.
    \item The DF method with hyperbolic ordering performs similarly to or outperforms the L-curve, NCP and W-GCV methods in all examples without exception.
    \item The DF method with elliptic ordering failed to produce acceptable solutions for the \texttt{Text2} problem with \(\alpha=10^{-4}\). Contrary to the other examples where the distance function \ref{eq:DistNormGKB} with this ordering has a minimum, in this example it has neither a minimum nor even an inflection point and therefore the right stopping iteration could not be found with this ordering.
    \item The optimal MSD values for the PLS problem \ref{eq:GKBprob} and the projected Tikhonov problem \ref{eq:TikhMinProb} are almost identical in all examples, as expected from the discussion in \ref{sec:TikhRegProj}.
\end{enumerate}
 Overall, we can conclude that the DF criterion with hyperbolic ordering for estimation of optimal stopping iteration is accurate, robust and outperforms state-of-the-art methods. 


\bibliographystyle{plain}
\bibliography{library}
\end{document}